\documentclass{article}
\usepackage[letterpaper, portrait, margin=2.5cm]{geometry}

\usepackage{verbatim}
\usepackage[utf8]{inputenc}
\usepackage{amsmath,amssymb,amsthm}
\usepackage{tikz-cd}
\usepackage{bbm}
\usepackage{hyperref}
\usepackage{comment}
\usepackage{stmaryrd}
\usepackage[new]{old-arrows}
\usepackage{graphicx}
\graphicspath{ {./Images/} }
\usepackage{todonotes}

\usepackage [english]{babel}
\usepackage [autostyle, english = american]{csquotes}
\MakeOuterQuote{"}

\usepackage[maxbibnames=99]{biblatex}
\newcommand{\PP}{\mathbb P}

\newcommand{\AAA}{\mathbb A}

\newcommand{\mce}{\mathcal E}

\newcommand{\mcl}{\mathcal L}

\newcommand{\lbm}{\left[ \begin{matrix}}
\newcommand{\rem}{\end{matrix} \right]}
\newcommand{\SL}{\sum\limits}

\newcommand{\VV}{\mathbb{V}}

\theoremstyle{definition}

\newtheorem{remark}{Remark}

\title{On the Application of Tschirnhaus Transformations to the Reduction of Algebraic Equations}
\author{A Translation by Alexander J. Sutherland}
\date{\today}

\begin{document}

\maketitle

\vspace{16pt}

\section{Original Bibliographic Information}
\subsection{German}
\begin{itemize}
	\item \textbf{Author:} Anders Wiman
	\item \textbf{Title:} \"Uber die Anwendung der Tschirnhausen-Transformation auf die Reduktion algebraischer Gleichungen
	\item \textbf{Year:} 1927
	\item \textbf{Language:} German
	\item \textbf{Publisher:} Nova Acta Regiae Societatis Scientiarum Upsaliensis (Nova Acta R. Soc. scient. Uppsala)
	\item \textbf{Note:} Der K\"onigi. Societ\"at der Wissenschaften Zu Uppsala Mitgeteilt
\end{itemize}

\subsection{English}
\begin{itemize}
	\item \textbf{Author:} Anders Wiman
	\item \textbf{Title:} On the Application of Tschirnhaus Transformations to the Reduction of Algebraic Equations
	\item \textbf{Year:} 1927
	\item \textbf{Language:} German
	\item \textbf{Publisher:} New Proceedings of the Royal Society of Scientists of Uppsala
	\item \textbf{Note:} Notice to the Royal Society of Scientists of Uppsala on 06 May 1927
\end{itemize}

\vfill

\begin{center}
This work was supported by the National Science Foundation under Grant No. DMS-1944862.
\end{center}

\newpage
\section{The Translation}

\subsection*{Part 1}
\footnote{Translator's Note: This is a translation of the original mathematics. In particular, errors in the text have not been fixed. The errors in question come from considering intersections in affine spaces instead of in projective spaces. Throughout this translation, there are additional footnotes with the identifier ``Translator's Footnote:.'' These footnotes refer to remarks in Section 3 in which the translator provides additional mathematical commentary.}We consider a general equation of $n^{th}$ degree:
\begin{equation}
	x^n + c_1x^{n-1} + \cdots + c_n = 0
\end{equation}

\noindent
with roots $x_1,\dotsc,x_n$. We then apply a Tschirnhaus transformation, which has the general form
\begin{equation}
	y = a_0 + a_1x + \cdots + a_{n-1}x^{n-1}.
\end{equation}

\noindent
This will transform (1) into an equation
\begin{equation}
	y^n + C_1y^{n-1} + \cdots + C_n = 0,
\end{equation}

\noindent
where the coefficients $C_i$ are all homogeneous functions of degree $i$ in the variables $a_{\nu}$ and include all such functions of weight up to $i$ in the coefficents $c_i$. Tschirnhaus hoped to use this type of transformation to convert equation (1) to the binomial form in such a way that the determination of parameters $a_{\nu}$ should require the solution of equations of degree at most $n-1$. As is well known, this is not the case, even though there have many attempts for the degree 5; this will never be the case, just like for [the problem of] trisecting an angle. However, the situation is completely different if the problem is formulated in the following way: Is it possible to satisfy the conditions
\begin{equation}
	C_i(a_0,a_1,\dotsc,a_{n-1})=0, \quad i=1,\dotsc,m
\end{equation}

\noindent
in such a way that determining the parameters $a_{\nu}$ only requires equations of degree up to $m$, \textit{when $n$ is sufficiently large?} As a result of the following treatment for the case $m=4$, it should not appear doubtful that question should also be decided in the affirmative for larger $m$. However, the general problem of determining the lower bound on $n$ associated to each $m$ appears to be very complex.

\subsection*{Section 2}
Observe that
\begin{align*}
	C_1(a_0,\dotsc,a_{n-1}) = na_0 + \cdots.
\end{align*}

If $C_1=0$, then $a_0$ is expressed linearly in the other parameters. The coefficients $C_i, (i=2,\dotsc,n-1)$ are then homogeneous functions of degree $i$ in the parameters $a_1,\dotsc,a_{n-1}$. In order to [find a point that will] satisfy a single condition
\begin{equation*}
	C_x = 0, \quad (x>1)
\end{equation*}

\noindent
it is evident that it is only necessary to find an intersection of the hypersurface $C_x=0$ with an arbitrary straight line to solve an equation of degree $x$. If all the roots of (1) are real, then you cannot get a real solution for $x=2$ because the hypersurface
\begin{equation}
	\SL_{i=1}^n y_i^2 = C_2 = 0
\end{equation}

\noindent
has only the trivial 0. In contrast, there are always real points on the surface $C_3=0$. Indeed, as you can easily see, the same is the case for all surfaces $C_x=0 \ (x>2)$ if $x$ is an even number. \\

If one has $n \geq 5$, then one obtains the solution of (4) for $m=3$ by the well-known Bring-Jerrard transformation, which is illustrated geometrically by F. Klein \footnote{We refer to the in-depth treatment of F. Klein, \textit{Lectures on the Icosahedron and the Solution of the Equation of Fifth Degree}, Leipzig, 1884.} in the following way.\footnote{Translator's footnote: See Remark \ref{rem:MainArgumentofSection2} for the translator's summary of this argument.} First, a point $P$ is obtained on the surface $C_2=0$, which, as noted above, can be done using a square root. If $n>5$, we then consider a three-dimensional space $R_3$ which is tangent [to $C_2=0$] at $P$, which then has a hypersuface of degree 2 in common with $C_2=0$. A second square root is now required in order to select one of the two generators of this surface going through $P$. Determining an intersection of one of these generators with $F_3=0$ requires the solution of a degree 3 equation. Although only one pair of imaginary roots need to occur in the real equations $C_1 = C_2 = C_3 = 0$, at least two pairs of imaginary roots must exist to actually execute this transformation. Otherwise, there is no real line at $C_2=0$.  This is due to the fact that if one converts
\begin{align*}
	C_2(a_1,\dotsc,a_{n-1}) = \SL_{i=1}^n y_i^2
\end{align*}

\noindent
to a sum of $n-1$ real squares, one gets the $\lambda$ with the sign $-$, where $2\lambda$ denotes the number of imaginary roots.

\subsection*{Section 3}
Now, let $n>5$. We assume that a point $P$ that lies on both $C_2=0$ and $C_3=0$ has been determined by the procedure given above. The associated coordinates are $a_1^{(0)},\dotsc,a_{n-1}^{(0)}$. We write
\begin{align*}
	\alpha_i - \alpha_i^{(0)} = \beta_i, \quad i=1,\dotsc,n-2
\end{align*}

\noindent
and reconstruct $C_2$ and $C_3$ in terms of $\beta_1,\dotsc,\beta_{n-2}$. \footnote{One can assume that $a_{n-1}^{(0)} \not= 0$, after possibly changing the indices, and then write $a_{n-1}=1$} In this manner, by summing the terms with same total degree in the $\beta_i$, we get:
\begin{align}
	C_2 &= \phi_1 + \phi_2\\
	C_3 &= \psi_1 + \psi_2 + \psi_3 \nonumber
\end{align}

We want to solve the present problem in such a way that we determine a straight line going through $P$ that lies on both the surfaces $C_2=0$ and $C_3=0$. If this is successful, then the further condition $C_4=0$ only requires the solution of an equation of the fourth degree. We denote the space with homogeneous coordinates $a_1,\dotsc,a_{n-1}$ as a $R_{n-2}$, in accordance with the number that is its dimension. \\

The first conditions to be introduced are
\begin{equation}
	\phi_1 = \psi_1 = 0.
\end{equation}

\noindent
The space $R_{n-2}$ is then reduced to a $R_{n-4}$.\footnote{Translator's footnote: See Remark \ref{rem:MainArgumentofSection3} for the translator's summary of this argument.} If we consider the straight lines through $P$ as elements of the space, the point space $R_{n-4}$ has only $n-5$ dimensions. From this perspective, we refer to it as $\mcl_{n-5}$ and consider the subvarieties $\phi_2=0$, $\psi_2=0$, and $\psi_3=0$ inside it. If one takes any plane in this line space $\mcl_{n-5}$, then $\phi_2=0$ and $\psi_2=0$ have four elements in common - that is, four straight line generators. If $n>7$, then a common straight line of $\phi_2=0$ and $\psi_2=0$ can be determined by solving a fourth degree equation. We denote such a straight line by $\ell_1$. \\

If we can complete the task in such a way that we get a plane tangent to $\phi_2=0$ and $\psi_2=0$, then this result will be solved, as this plane intersects $\psi_3=0$ in three straight lines, so that everything else comes down to the solution of a degree three equation. The equations of the hyperplanes, which meet the subvarieties $\phi_2=0$ and $\psi_2=0$ along the straight line $\ell_1$ are given by
\begin{equation}
 \phi_1^{(1)} = 0, \quad \psi_1^{(1)} = 0. 
\end{equation}

\noindent
We assume the relations (8) are satisfied. The $R_{n-4}$ discussed above then reduces to an $R_{n-6}$. If we restrict to this $\PP^{n-6}$, then elements of $\phi_2=0$ and $\psi_2=0$ appear as planes through $\ell_1$. However, there are exceptions for when $n=7$, as $R_{n-6}$ coincides with $\ell_1$, and for when $n=8$, as $\phi_2=0$ and $\psi_2=0$ are reduced to the double-counted straight line $\ell_1$. We now replace the point space $R_{n-6}$ with a space whose elements are the planes that contain the line $\ell_1$. This plane space obviously has dimension $n-8$ and is denoted by $\mce_{n-8}$. We now have the intersection of the quadratic [algebraic] manifolds $\phi_2=0$ and $\psi_2=0$ in $R_{n-6}$ in the plane space $\mce_{n-8}$. If $n \geq 10$, one concludes that only a fourth degree equation has to be solved to determine a common plane of $\phi_2=0$ and $\psi_2=0$ in $\mce_{n-8}$. Therefore, we have the theorem:\\

\textit{If $n \geq 10$, the general equation (1) can be reduced to the form}
\begin{equation}
	y^n + C_5y^{n-5} + \cdots + C_n = 0
\end{equation}

\noindent
\textit{by a transformation (2) in such a way that determining the parameters $a_i$ requires only the solution of a finitely many quartic, cubic, and quadratic irrationalities.} \\

Although an equation of the form (9) can be achieved with only two pairs of imaginary roots, at least three pairs of imaginary roots of (1) must exist for this transformation to be completed using only real numbers. It is only under this condition that the hypersurface $C_2=0$ contains real planes. However, if one wishes for not only necessary, but also sufficient conditions here, this does seem to be possible without fairly in-depth discussions.

\subsection*{Section 4}
The case $n=9$ was examined by D. Hilbert in a recently published work.\footnote{\"Uber die Gleichungen neuten Grades, Math. Ann. 97, S. 243 (1926)} It is first demonstrated how one can determine a three-dimensional space $R_3$ which is completely contained in the hypersurface $C_2=0$ and then one considers the degree three surface $F_3$ in this $R_3$ which is cut out by the hypersurface $C_3=0$. As is known, this surface $F_3$ only depends on four fundamental parameters. This is also particularly evident when one transforms only the left term of the equation to a sum of five cubes, for which it is necessary to solve a fifth degree equation. Hilbert now puts the general equation of ninth degree in the form
\begin{equation}
	y^9 + C_5y^4 + C_6y^3 + C_7y^2 + C_8y + C_9 = 0
\end{equation}

\noindent
by first determining one of the 27 straight lines on the surface $F_3$ and then intersecting one of these straight lines with the hypersurface $C_4=0$. Both the equation (10), where one can easily set $C_9=1$, as well as the equation of degree 27 are functions of only four parameters. The result of this is that the solution of the general equation of ninth degree only requires algebraic functions of four arguments in such a way that "one can get by with functions of one argument, sums, and two special functions of four arguments". \\

It can be shown that, if the general equation of ninth degree is reduced to the form (10), then there is no need to solve an equation of degree higher than five, so that one of the special functions of four arguments above is unnecessary. At the end, we generalize our task set at (4) by using auxiliary equations larger than $m$, but with the restriction that each degree will still always be $<n$. \\

As in the previous case, we determine a point $P$ on both $C_2=0$ and $C_3=0$ and still suppose the conditions (8) hold. Since $n=9$ here, the left terms of the subcone $\phi_2=0$ and $\psi_2=0$ can be written in five homogeneous coordinates, such as $z_1,z_2,z_3,z_4,z_5$. We now look for the self-conjugate pentahedron common to both $\phi_2$ and $\psi_2$, which corresponds to solving a fifth degree equation.\footnote{Translator's footnote: See Remark \ref{rem:PencilsOfQuadrics} for more exposition.} If this pentahedron is assumed to be a coordinate pentahedron, we have expressions for $\phi_2$ and $\psi_2$ of the form:
\begin{align}
	\phi_2 &= a_1z_1^2 + a_2z_2^2 + a_3z_3^2 + a_4z_4^2 + a_5z_5^2\\
	\psi_2 &= b_1z_1^2 + b_2z_2^2 + b_3z_3^2 + b_4z_4^2 + b_5z_5^2. \nonumber
\end{align}

One can now eliminate any of the five variables from $\phi_2=0$ and $\psi_2=0$ and thus obtain five relations, one of which we write in the form
\begin{equation}
	c_1z_1^2 + c_2z_2^2 + c_3z_3^2 + c_4z_4^2 = 0. 
\end{equation}

\noindent
We take (12) as the equation of a degree two hypersurface and try to determine the corresponding straight line generators, which only requires taking square roots. The straight lines in this generating set are assigned to the values of a a parameter $\lambda$ and, likewise, the points of a specified generator are assigned to the values of another parameter $t_1$ and the their coordinates $z_1,z_2,z_3,z_4$ can be expressed linearly in both $\lambda$ and $t$. According to (11), we get the relation for $z_5$:
\begin{equation}
	z_5^2 = a_2(\lambda)t^2 + b_2(\lambda)t + c_2(\lambda).
\end{equation}

Since the elements of (11) are straight lines through $P$, it can be seen that a generator of (12) corresponds to a two-dimensional cone whose apex is $P$. Since this cone must have six generators in common with the hypersurface $\psi_3=0$, it follows that one need not use auxiliary equations of degree more than six when reducing the general equation of ninth degree to the form (10). 

\subsection*{Section 5}
This matter can be simplified by looking for a value of $\lambda$ such that the two-dimensional cone splits into two planes. We need only solve the degree four equation
\begin{equation}
	\left[ b_2(\lambda) \right]^2 - 4a_2(\lambda)c_2(\lambda) = 0.
\end{equation}

\noindent
The generating system of (12) has four other conjugate pairs of planes, so that the total number of planes common to the subcones $\phi_2=0$ and $\psi_2=0$ is 16. According to the five different relations in four variables of the form (12), each of these planes can be paired with five others. Two planes that can be paired together intersect each other in a straight line. If not, they only intersect at $P$. \\

The theory of the intersection of the two cones $\phi_2=0$ and $\psi_2=0$ is indeed well known.\footnote{Translator's footnote: See Remark \ref{rem:Degree4DelPezzoSurfaces} for the translator's summary of this argument.} They system of equations (11) is often used to study properties of a degree four surface with a double conic section. The 16 common planes of the two subcones correspond to the 16 straight lines lying on such a surface.\\

We have now demonstrated \textit{that the general equation of ninth degree can also be converted to the form (9) without having to presuppose the solution of auxiliary equations which each depend on more than one parameter.} Among the auxiliary equations, however, there is one of fifth-degree, so we cannot get by with square roots and cube roots, as we can for $n>9$. However, it seems to be difficult to prove strictly that the latter is not possible at all for $n=9$.

\newpage
\section{Translator's Notes}
\begin{remark}\label{rem:MainArgumentofSection2} \textbf{(Main Argument of Section 2)}\\
Let $\AAA^n$ be the affine space of Tschirnhaus transformations and $\PP^{n-1}$ its projectivization. $\VV(C_1)$ is a hyperplane and thus $\VV(C_1) \cong \PP^{n-2}$. A rational point $P$ of $\VV(C_1) \cap \VV(C_2)$ can be determined over a quadratic extension of the base field. Then, the tangent hyperplane $T$ at $P$ can be computed rationally and will have dimension at least 4 in $\VV(C_1) = \PP^{n-2}$ if $n>6$. Hence $\VV(C_1) \cap \VV(C_2) \cap T$ is a quadric that is singular at $P$ in $\VV(C_1) \cap T \cong \PP^{n-3}$. Consequently, a line in this cone can be determined by solving a quadratic polynomial and it suffices to intersect this line with $F_3=0$. 
\end{remark}

\begin{remark}\label{rem:MainArgumentofSection3} \textbf{(Main Argument of Section 3)}\\
To re-state Wiman's approach, this $\PP^{n-4}$ is obtained by projectivizing and then considering $\VV(C_1) \cap \VV(\phi_1) \cap \VV(\psi_1)$ inside $\PP^{n-1}$. By shifting $P$ to the origin (e.g. $[0:\cdots:0:1]$), Wiman then uses a classical correspondence to consider the pencil of lines through $P$ and identifies it as $\mcl_{n-5} \cong \PP^{n-5}$. Note that $\phi_2, \psi_2,$ and $\psi_3$ induce hypersurfaces of the same degree in $\mcl_{n-5}$. If $n \geq 7$, a point $Q$ of $\VV(\phi_2) \cap \VV(\psi_2) \subseteq \mcl_{n-5}$ can be determined by solving a quartic equation. Moreover, by construction, the line $\ell_1$ determined by $P$ and $Q$ lies in $\VV(C_1) \cap \VV(C_2)$ in the ambient space. \\

Now, consider the tangent hyperplanes of $\VV(\phi_2), \VV(\psi_2) \subseteq \PP^{n-4}$ defined by the polynomials $\phi_1^{(1)}$ and $\psi_1^{(1)}$. Consider the $\PP^{n-6}$ given by $\VV(C_1) \cap \VV(\phi_1) \cap \VV(\psi_1) \cap \VV(\phi_1^{(1)}) \cap \VV(\psi_1^{(1)})$. Every point not on $\ell_1$ in $\mce_{n-8} = \VV(\phi_2) \cap \VV(\psi_2) \subseteq \PP^{n-6}$ determines a plane on $\VV(\phi_2) \cap \VV(\psi_2)$ in the ambient space; this is possible when $n \geq 9$. Determining such a point $Q'$ in $E_{n-8}$ determines a line $\ell_2 \subseteq \VV(\phi_2) \cap \VV(\psi_2) \subseteq \mcl_{n-5}$ and thus a point $Q'$ of $\ell_2 \cap \VV(\psi_3) \subseteq \mcl_{n-5}$ can be determined by solving a cubic equation. However, the line determined by $P$ and $Q'$ lies inside $\VV(C_1) \cap \VV(C_2) \cap \VV(C_3)$ in the ambient space and thus a point of $\VV(C_1) \cap \VV(C_2) \cap \VV(C_3) \cap \VV(C_4)$ can be determined by solving a quartic equation.
\end{remark}

\begin{remark}\label{rem:PencilsOfQuadrics} \textbf{(Pencils of Quadrics)}\\
The forms $\phi_2$ and $\psi_2$ define a pencil of quadratic forms in the five variables $z_1,\dotsc,z_5$. The singular fibers of the pencil are given by the roots of the discriminant, which is a polynomial of degree 5. Hence, determining a singular quadric in the pencil corresponds to solving a degree 5 polynomial. Wiman then again uses the observation that a singular quadric is a cone.
\end{remark}

\begin{remark}\label{rem:Degree4DelPezzoSurfaces} \textbf{(Degree 4 del Pezzo Surfaces)}\\
Here Wiman is observing the fact that the intersection of two quadrics in $\PP^4$ is a degree 4 del Pezzo surface.
\end{remark}


\end{document}